\documentclass[a4paper,12pt]{article}
\usepackage{amssymb} \usepackage{amsfonts}
\usepackage{amsmath}
\usepackage[utf8x]{inputenc}
\usepackage{tikz,relsize,mathtools}
\usepackage[cmtip,arrow]{xy}
\usepackage{pb-diagram,pb-xy}
\newcommand{\F}{{\mathbb F}}
\newcommand{\Z}{{\mathbb Z}}
\newcommand{\Q}{{\mathbb Q}}
\newcommand{\R}{{\mathbb R}}
\newcommand{\C}{{\mathbb C}}

\newcommand{\N}{{\mathbb N}}
\newcommand{\M}{{\mathbb S}}
\begin{document}
\title{Flat manifolds with holonomy group $\Z_{2}^{k}$ of diagonal type}
\author{A. G\c{a}sior, A. Szczepa\'nski}
\date{\today}
\maketitle
\section{Introduction}
Let $M^n$ be a flat manifold of dimension $n.$
By definition, this is a compact connected, Riemannian manifold without boundary
with sectional curvature equal to zero. From the theorems of Bieberbach (\cite{Ch})
the fundamental group
$\pi_{1}(M^{n}) = \Gamma$ determines a short exact sequence:
\begin{equation}\label{ses}
0 \rightarrow \Z^{n} \rightarrow \Gamma \stackrel{p}\rightarrow
G \rightarrow 0,
\end{equation}
where
$\Z^{n}$ is a torsion free abelian group of rank $n$ and
$G$ is a finite group which
is isomorphic to the holonomy group of $M^{n}.$
The universal covering of $M^{n}$ is the Euclidean space $\R^{n}$
and hence $\Gamma$
is isomorphic to a discrete cocompact subgroup
of the isometry group $Isom(\R^{n}) = O(n)\times\R^n = E(n).$
Conversely, given a short exact sequence of the form (\ref{ses}), it is known that
the group $\Gamma$ is (isomorphic to) the fundamental group of a flat manifold if and only if
$\Gamma$ is torsion free.
In this case $\Gamma$ is called a Bieberbach group.
We can define a holonomy representation $\phi:G\to GL(n,\Z)$ by the formula:
\begin{equation}\label{holonomyrep}
\forall g\in G,\phi(g)(e_i) = \tilde{g}e_i(\tilde{g})^{-1},
\end{equation}
where $e_i\in\Gamma$ are generators of $\Z^n$ for $i=1,2,...,n,$ and $\tilde{g}\in\Gamma$
such that $p(\tilde{g})=g.$ In this article we shall consider only the case
\begin{equation}\label{diagonal}
G = \Z_{2}^{k}, 1\leq k\leq n-1,
\hskip 2mm\text{with}\hskip 2mm \phi(\Z_{2}^{k})\subset D\subset GL(n,\Z),
\end{equation}
where $D$ is the group of all diagonal matrices.
We want to consider relations between two families of flat manifolds with the above property (\ref{diagonal}):
the family ${\cal RBM}$ of real Bott manifolds and the family ${\cal GHW}$ of generalized 
Hantzsche-Wendt manifolds.
In particular, we shall prove (Proposition \ref{crossing}) that the intersection ${\cal GHW}\cap {\cal RBM}$ is not empty.
\vskip 1mm
\noindent
In the next section we consider some class of real Bott manifolds without $\operatorname{Spin}$ and $\operatorname{Spin}^{\C}$ structure.
There are given conditions (Theorem~\ref{main}) for the existence of such structures.
As an application a list of all 5-dimensional oriented real Bott manifolds without $\operatorname{Spin}$ structure is given
, see Example~\ref{ex2}.
In this case we generalize the results of L. Auslaneder and R. H. Szczarba, \cite{AS} from 1962, cf. Remark~\ref{remark}.
At the end we formulate a question about cohomological rigidity of ${\cal GHW}$ manifolds.
\section{Families}
\subsection{Generalized Hantzsche-Wendt manifolds}
We start with the definition of generalized Hantzsche-Wendt manifold.
\newtheorem{defi}{Definition}
\begin{defi} {\em (\cite[Definition]{RS})}
A generalized Hantzsche-Wendt manifold (for short ${\cal GHW}$-manifold) 
is a flat manifold of dimension $n$ with holonomy group
$(\Z_2)^{n-1}.$
\end{defi}
Let $M^n\in {\cal GHW}.$ In \cite[Theorem 3.1]{RS} it is proved that
the holonomy representation (\ref{holonomyrep}) of $\pi_1(M^n)$ satisfies (\ref{diagonal}).
\vskip 1mm
\noindent
The simple and unique example of an oriented $3$-dimensional generalized Hantzsche-Wendt manifold
is a flat manifold which was considered for the first time by W. Hantzsche and H. Wendt in 1934, \cite{HW}.
\vskip 1mm
\noindent
Let $M^n\in {\cal GHW}$ be an oriented, $n$-dimensional manifold (a HW-manifold for short).
In 1982, see \cite{RS}, the second author proved that for odd $n\geq 3$ and for all $i, H^i(M^n,\Q) \simeq H^i(\M^n,\Q),$
where $\Q$ are the rational numbers and $\M^n$ denotes the $n$-dimensional sphere. Moreover,
for $n\geq 5$ the
commutator subgroup of the fundamental group $\pi_{1}(M^n) = \Gamma$ is equal to the
translation subgroup ($[\Gamma,\Gamma] = \Gamma\cap\R^n$), \cite{P}.
The number $\Phi(n)$ of affine non equivalent HW-manifolds of
dimension $n$ growths exponetially, see \cite[Theorem 2.8]{MR}, and
for $m\geq 7$ there exist many isospectral manifolds 
non pairwise homeomorphic, \cite[Corollary 3.6]{MR}. The
manifolds have an interesting connection with Fibonacci groups
\cite{S1} and the theory of quadratic forms over the field $\F_2,$
\cite{S2}. HW-manifolds have no $\operatorname{Spin}$-structure, \cite[Example 4.6 on page 4593]{MP}. 
\vskip 1mm
The (co)homology groups and cohomology rings with coefficients in $\Z$ or $\Z_2,$ of generalized Hantzsche-Wendt manifolds
are still not known, see \cite{CMR} and \cite{DP}. 
%Finally, let us mention about
%properties related to the theory of fixed points. HW-manifolds
%satisfy so called Anosov relation. This means for any continious map
%$f:M^n\to M^n, \mid L(f)\mid = N(f),$ where $L(f)$ is the Lefschetz
%number of $f$ and $N(f)$ is the Nielsen number of $f,$ see
%\cite{DR}.
We finish this overview with an example of generalized Hantzsche-Wendt manifolds which have been known
already in 1974. 
\newtheorem{ex}{Example}
\begin{ex}
Let $M^n = \R^n/\Gamma_n, n\geq 2$ be manifolds defined in
{\em  \cite{LS} (}see also {\em \cite[page 1059]{RS})},
where $\Gamma_n\subset E(n)$ is generated by
$\gamma_{0} = (I=\text{id},(1,0,...,0))$ and
\begin{equation}\label{gener1}
\gamma_{i} = \left(\left[
\begin{matrix}
1&0&0&.&.&...&0\\
0&1&0&.&.&...&0\\
.&.&.&.&.&...&\\
0&...&1&0&0&...&0\\
0&...&0&-1&0&...&0\\
0&...&0&0&1&...&0\\
.&.&.&.&.&...&\\
0&.&...&0&0&0&1
\end{matrix}\right], \begin{pmatrix}
0\\
0\\
...\\
0\\
0\\
\frac{1}{2}\\
...\\
0
\end{pmatrix}\right)\in E(n),
\end{equation}
where the $-1$ is placed in the $(i,i)$ entry and the $\frac{1}{2}$ as an $(i+1)$ entry, $i = 1,2,...,n-1.$
$\Gamma_2$ is the fundamental group of the Klein bottle.
\end{ex}

\subsection{Real Bott manifolds}
We follow \cite{CMO}, \cite{KM} and \cite{Nazra}.
To define the second family let us introduce a sequence of $\R P^1$-bundles
\begin{equation}\label{tower}
M_{n}\stackrel{\R P^1}\to M_{n-1}\stackrel{\R P^1}\to...\stackrel{\R P^1}\to M_{1}\stackrel{\R P^1}\to M_0 =
\{ \text{a}\hskip 1mm \text{point}\}
\end{equation}
such that $M_{i}\to M_{i-1}$ for $i = 1,2,...,n$ is the projective bundle of a Whitney sum
of a real line bundle $L_{i-1}$ and the trivial line bundle over $M_{i-1}.$ 
We call the sequence (\ref{tower}) a {\em real Bott tower} of height $n,$ \cite{CMO}.
\begin{defi}{\em (\cite{KM})}
The top manifold $M_n$ of a real Bott tower {\em (\ref{tower})} is called a real Bott manifold.
\end{defi}
Let $\gamma_i$ be the canonical line bundle over $M_i$ and set $x_i = w_1(\gamma_j).$
Since $H^1(M_{i-1},\Z_2)$ is additively generated by $x_1,x_2,..,x_{i-1}$ and $L_{i-1}$
is a line bundle over $M_{i-1},$ one can uniquely write
\begin{equation}\label{w1}
w_1(L_{i-1}) = \Sigma_{k=1}^{i-1} a_{k,i}x_k
\end{equation}
with $a_{k,i}\in \Z_2 = \{0,1\}$ and $i = 2,3,...,n.$

From above $A = [a_{ki}]$ is an upper triangular matrix
\footnote{$a_{k,i} = 0$ unless $k < i.$}
of size $n$ 
whose diagonal entries are $0$ and other entries are either
$0$ or $1.$ 
Summing up, we can say that the tower (\ref{tower}) is completly determined by the matrix $A.$
\vskip 1mm
\noindent
From \cite[Lemma 3.1]{KM} we can consider any real Bott manifold $M(A)$ in
the following way.  
Let $M(A) = \R^n/\Gamma(A),$
where $\Gamma(A)\subset E(n)$ is generated by elements
\begin{equation}\label{gener}
s_{i} = \left(\left[
\begin{matrix}
1&0&0&.&.&...&0\\
0&1&0&.&.&...&0\\
.&.&.&.&.&...&\\
0&...&0&1&0&...&0\\
0&...&0&0&(-1)^{a_{i,i+1}}&...&0\\
.&.&.&.&.&...&\\
0&...&0&0&0&...&(-1)^{a_{i,n}}
\end{matrix}\right], \begin{pmatrix}
0\\
.\\
0\\
\frac{1}{2}\\
0\\
.\\
0\\
0
\end{pmatrix}\right)\in E(n),
\end{equation}
where $(-1)^{a_{i,i+1}}$ is placed in $(i+1, i+1)$ entry and $\frac{1}{2}$ as an $(i)$ entry, $i = 1,2,...,n-1.$
$s_{n} = (I,(0,0,...,0,1))\in E(n).$ 
From \cite[Lemma 3.2,3.3]{KM} $s_{1}^{2},s_{2}^{2},...,s_{n}^{2}$ commute with each
other and generate a free abelian subgroup $\Z^n.$ It is easy to see that it is not always a maximal abelian subgroup
of $\Gamma(A).$ Moreover, we have the following commutative diagram
$$\begin{diagram}
\node{0}\arrow{e}\node{N}\arrow{e}\node{\Gamma(A)}\arrow{s,=}\arrow{e}\node{\Z_{2}^{k}}\arrow{e}\node{ 0}\\
\node{0}\arrow{e}\node{\Z^n}\arrow{n,r}{i}\arrow{e}\node{\Gamma(A)}\arrow{e}\node{\Z_{2}^{n}}\arrow{n,r}{p}\arrow{e}\node{0}
\end{diagram}$$
where $k = rk_{\Z_2}(A), N$ is the maximal abelian subgroup of $\Gamma(A)$,
and $p:\Gamma(A)/\Z^n\to\Gamma(A)/N$ is a surjection induced by the inclusion $i:\Z^n\to N.$
From the first Bieberbach theorem, see \cite{Ch}, $N$ is a subgroup of all translations of $\Gamma(A)$ i.e.
$N = \Gamma(A)\cap\R^n = \Gamma(A)\cap \{(I,a)\in E(n)\mid a\in\R^n\}.$
\vskip 1mm
\begin{defi} {\em (\cite{CMO})}
A binary square matrix $A$ is a Bott matrix if $A = PBP^{-1}$ for a permutation
matrix $P$ and a strictly upper triangular binary matrix $B.$
\end{defi}
Let ${\cal B}(n)$ be the set of Bott matrices of size $n.$
\footnote{Sometimes ${\cal B}(n)$ is defined to be the set of strictly upper triangular binary matrices of size $n.$}
Since two different upper triangular matrices  $A$ and $B$ may produce (affinely)
diffeomorphic ($\sim$) real Bott manifolds $M(A), M(B),$ see \cite{CMO} and \cite{KM}, there are 
three operations on ${\cal B}(n),$ denoted by (Op1), (Op2) and (Op3),
such that $M(A)\sim M(B)$ if and only if the matrix $A$ can be transformed into $B$ 
through a sequence of the above operations, see \cite[part 3]{CMO}.
The operation (Op1) is a conjugation by a permutation matrix,
\vskip 1mm
\noindent
(Op2) is a bijection $\Phi_k:{\cal B}(n)\to {\cal B}(n)$
\begin{equation}\label{op2}
\Phi_k(A)_{\ast,j} := A_{\ast,j} + a_{kj}A_{\ast,k},
\end{equation}
for $k,j\in\{1,2,...,n\}$ such that $\Phi_k\circ\Phi_k$ = $1_{{\cal B}(n)}.$
\vskip 1mm
\noindent
Finally (Op3) is, for distinct $l,m\in\{1,2,..,n\}$ and the matrix $A$ with $A_{\ast,l} = A_{\ast,m}$ 
\begin{equation}\label{op3}
\Phi^{l,m}(A)_{i,\ast} := \left\{ \begin{array}{ll}                                                                             
A_{l,\ast} + A_{m,\ast} & \text{if}\,\, i=m\\ A_{i,\ast} & \text{otherwise}
\end{array}\right.
\end{equation}
\noindent
Here $A_{\ast,j}$ denotes $j$-th column and $A_{i,\ast}$ denotes $i$-th row of the matrix $A.$
\newtheorem{theo}{Theorem}
\vskip 5mm
Let us start to consider the relations between these two classes of flat manifolds.
We start with an easy observation
$${\cal RBM}(n)\cap {\cal GHW}(n) = \{M(A)\mid rank_{\Z_2}A = n-1\} =$$
$$= \{M(A)\mid a_{1,2}a_{2,3}...a_{n-1,n} = 1\}.$$
\noindent
These manifolds are classified in \cite[Example 3.2]{CMO} and for $n\geq 2$
\begin{equation}\label{rbhw}
\#({\cal RBM}(n)\cap {\cal GHW}(n)) = 2^{(n-2)(n-3)/2}.
\end{equation}
\vskip 2mm
\noindent
There exists the classification, see \cite{RS} and \cite{CMO}, 
of diffeomorphism classes of ${\cal GHW}$ and ${\cal RBM}$ manifolds in low dimensions.
For dim $\leq 6$ we have the following table.
\vskip 5mm
\begin{tabular}{|c||c|c||c|c||c|} 
\hline
%lub np. {lcrrr}, {|rr|} itd. w zależnoÅ od liczby żÄdanych kolumn i ich justowania
&\multicolumn{2}{c||}{number of}&\multicolumn{2}{c||}{number of}&number of \\
dim&\multicolumn{2}{c||}{$GHW$ manifolds}&\multicolumn{2}{c||}{$RBM$ manifolds}&$GHW\cap RBM$ manifolds\\
\hline
&total&oriented&total&oriented&total\\
\hline
1&0&0&1&1&0\\
\hline
2&1&0&2&1&0\\
\hline
3&3&1&4&2&1\\
\hline
4&12&0&12&3&2\\
\hline
5&123&2&54&8&8\\
\hline
6&2536&0&472&29&$64$\\
\hline
\end{tabular}
\vskip 5mm
\newtheorem{prop}{Proposition}
\begin{prop}\label{crossing}
$\Gamma_{n}\in {\cal GHW}\cap {\cal RBM}.$
\end{prop}
\noindent
{\bf Proof:}
It is enough to see that the group $(G,0)\Gamma_n(G,0)^{-1} = \Gamma(A),$
where $G = [g_{ij}], 1\leq i,j\leq n,$ 
$$ 
g_{ij} := \left\{ \begin{array}{ll}
1 & \text{if}\,\, j = n-i+1\\ 0 & \text{otherwise}
\end{array}\right.$$
and $A = [a_{ij}], 1\leq i,j\leq n,$
with 
$$
a_{ij} := \left\{ \begin{array}{ll}
1 & \text{if}\,\, j = i+1\\ 0 & \text{otherwise}
\end{array}\right.$$
\vskip 5mm
\hskip 122mm $\Box$
 \section{Existence of $\operatorname{Spin}$ and $\operatorname{Spin}^{\C}$ structures on real Bott manifolds}
In this section we shall give some condition for the existence of $\operatorname{Spin}$ and $\operatorname{Spin}^{\C}$
structures on real Bott manifolds. 
We use notations from the previous sections.
There are a few ways to decide whether there exists a $\operatorname{Spin}$ structure 
on an oriented flat manifold $M^n,$ see \cite{GS}.
We start with the following. A closed oriented differential manifold $N$ has such a 
structure if and only if the second Stiefel-Whitney class $w_2(N) = 0.$
In the case of an oriented real Bott manifold $M(A)$ we have the formula for $w_2.$
\vskip 1mm
\noindent
Recall, see \cite{KM}, that for the Bott matrix $A$
\begin{equation}\label{ring}
H^{\ast}(M(A);\Z_2) = \Z_2[x_1,x_2,...,x_n]/(x_{j}^{2} = x_{j}\Sigma_{i=1}^{n}a_{i,j}x_i\mid j = 1,2,...,n)
\end{equation}
as graded rings.
Moreover,
from \cite[(3.1) on page 3]{LS} the $k$-th Stiefel-Whitney class
\begin{equation}
w_k(M(A)) = (B(p))^{\ast}\sigma_{k}(y_1,y_2,...,y_{n})\in H^{k}(M(A);\Z_2) ,
\end{equation}
where $\sigma_k$ is the $k$-th elementary symmetric function, $$p:\pi_{1}(M(A))\to G\subset O(n)$$ a holonomy representation, $B(p)$ is a map
induced by $p$ on the classification spaces and $y_i \stackrel{(\ref{w1})}= w_1(L_{i-1}).$ 
Hence,
\begin{equation}\label{sw1}
w_{2}(M(A)) = \Sigma_{1\leq i< j\leq n} y_{i}y_{j}\in H^{2}(M(A);\Z_2).
\end{equation}
\vskip 1mm
\noindent
There exists a general condition, see \cite[Theorem 3.3]{CMR}, for the 
calculation of the second Stiefel-Whitney for flat manifolds with $(\Z_2)^k$ holonomy of diagonal type
but we prefer the above explicit formula (\ref{sw1}). 
\footnote{We use it in Example~\ref{ex2}.}
Its advantage follows from the knowledge of the cohomology ring (\ref{ring}) of real Bott manifolds.
\vskip 5mm
An equivalent condition for the existence of a $\operatorname{Spin}$ structure is as follows.
An oriented flat manifold $M^n$ (a Bieberbach group $\pi_1(M^n)=\Gamma$) 
has a $\operatorname{Spin}$ structure if and only if there exists
a homomorfism $\epsilon:\Gamma\to \operatorname{Spin}(n)$ such that $\lambda_n\epsilon = p.$
Here $\lambda_n:\operatorname{Spin}(n)\to SO(n)$ is the covering map, see \cite{GS}.
We have a similar condition, under assumption $H^2(M^n,\R) = 0,$ for 
the existence of $\operatorname{Spin}^{\C}$ structure, \cite[Theorem 1]{GS}.
In this case $M^n$ (a Bieberbach group $\Gamma$) has a $\operatorname{Spin}^{\C}$ structure 
if an only if there exists a homomorphism
\begin{equation}\label{eqx} 
\bar{\epsilon}:\Gamma\to \operatorname{Spin}^{\C}(n) 
\end{equation}
such that $\bar{\lambda_n}\bar{\epsilon} = p.$
$\bar{\lambda_n}:\operatorname{Spin}^{\C}(n)\to SO(n)$ is the homomorphism induced by $\lambda_n,$ see \cite{GS}.
\newtheorem{rem}{Remark}
We have the following easy observation.
If there existe $H\subset\Gamma,$ a subgroup of finite index, such that the finite covering $\tilde{M^n}$
with $\pi_1(\tilde{M^n}) = H$ has no $\operatorname{Spin}$ ({\em $\operatorname{Spin}^{\C}$}) structure,
then $M^n$ has also no such structure.
\vskip 4mm
\noindent
We shall prove.
\begin{theo}\label{main}
Let $A$ be a matrix of an orientable real Bott manifold $M(A)$ of dimension $n.$ 
\vskip 1mm
\noindent
{\bf I.} Let $l\in\N$ be an odd number.
If there exist $1\leq i < j\leq n$ and rows $A_{i,\ast}, A_{j,\ast}$
such that
\begin{equation}\label{cond1}
\#\{m\mid a_{i,m} = a_{j,m} = 1\} = l
\end{equation} and
\begin{equation}\label{cond2}
a_{i,j} = 0,
\end{equation}
then $M(A)$ has no $\operatorname{Spin}$ structure.
\vskip 1mm
\noindent
Moreover, if 
\begin{equation}\label{cond3}
\#\{J\subset\{1,2,...,n\}\mid \# J = 2,\Sigma_{j\in J} A_{\ast,j} = 0\} = 0,
\end{equation}
then $M(A)$ has no $\operatorname{Spin}^{\C}$ structure.
\vskip 5mm
\noindent
{\bf II.}
If there exist $1\leq i<j\leq n$ and rows 
$$A_{i,\ast} = (0,...,0,a_{i,i_{1}},....,a_{i,i_{2k}},0,...,0),$$
$$A_{j,\ast} = (0,...,0,a_{j,i_{2k+1}},...,a_{j,i_{2k+2l}},0,...,0)$$
such that 
$a_{i,i_{1}} = a_{i,i_{2}} = ... = a_{i,i_{2k}} = 1, a_{i,m} = 0$ for $m\notin\{i_1,i_2,...,i_{2k}\}$
$a_{j,i_{2k+1}} = a_{j,i_{2k+2}} = ... = a_{j,i_{2k+2l}} = 1, a_{j,r} = 0$ for 
$r\notin\{i_{2k+1},i_{2k+2},...,i_{2k+2l}\}$ 
and $l,k$ odd
then $M(A)$ has no $\operatorname{Spin}$ structure.
\end{theo}
{\bf Proof:}
From \cite[Lemma 2.1]{KM} the manifold $M(A)$ is orientable if and only if for any $i = 1,2,..,n,$
$$\Sigma_{k = i+1}^{n}a_{i,k} = 0\hskip 2mm \text{mod}\hskip 1mm 2.$$
\noindent
Assume that $\epsilon : \pi_1(M(A))\to \operatorname{Spin}(n)$ defines a $\operatorname{Spin}$ structure on $M(A).$
Let $a_{i,i_1},a_{i,i_2},...,a_{i,i_{2m}},a_{j,j_1},a_{j,j_2},...,a_{j,j_{2p}} = 1$
and let $s_i, s_j$ be elements of $\pi_1(M(A))$ which define rows $i, j$ of $A,$ see (\ref{gener}).
Then $$\epsilon(s_i) = \pm e_{i_1}e_{i_2}...e_{i_{2m}},$$
$$\epsilon(s_j) = \pm e_{j_1}e_{j_2}...e_{j_{2p}}$$ and
$$\epsilon(s_i s_j) = \pm e_{k_1}e_{k_2}...e_{k_{2r}}.$$
From (\ref{cond1}) $2r = 2m+2p-2l.$
Moreover $\epsilon(s_{i}^{2}) = (-1)^{m}, \epsilon(s_{j}^{2}) = (-1)^{p}$ and
$\epsilon((s_i s_j)^2) = (-1)^{m+p-l} = (-1)^{m+p+l}.$
Since from (\ref{cond2}) (see also \cite[Lemma 3.2]{KM})
$s_i s_j = s_j s_i$ we have $\epsilon((s_i)^2)\epsilon((s_j)^2) = \epsilon((s_i s_j)^2).$
Hence $$(-1)^{m+p} = (-1)^{m+p+l}.$$ This is impossible since $l$ is an odd number and we have
a contradiction.
\vskip 2mm
\noindent
For the existence of the $\operatorname{Spin}^{\C}$ structure 
it is enough to observe that the condition (\ref{cond3}) is
equivalent to equation $H^2(M(A),\R) = 0,$ see \cite[formula (8.1)]{CMO}.
Hence, we can apply the formula (\ref{eqx}). 
Let us assume that $\bar{\epsilon}:\pi_1(M(A)\to \operatorname{Spin}^{\C}(n)$
defines a $\operatorname{Spin}^{\C}$ structure. Using the same arguments as 
above, see \cite[Proposition 1]{GS}, we obtain a contradiction. This finished the proof of {\bf I.}
\vskip 3mm
\noindent
For the proof {\bf II} let us observe that $s_{i}^{2} = (s_{i}s_{j})^2.$
Hence $(-1)^{k} = \epsilon((s_{i})^{2} = \epsilon((s_{i}s_{j})^2) = (-1)^{k+l} = 1.$
This is impossible.
\vskip 3mm
\hskip 122mm $\Box$
\vskip 5mm
\noindent
In the above theorem rows of number $i$ and $j$ correspond to generators $s_i,s_j$ which define
a finite index subgroup $H\subset\pi_1(M(A)).$
It is a Bieberbach group with holonomy group $\Z_2\oplus\Z_2.$ We proved that $H$ (if it exists) has
no $\operatorname{Spin}$ ($\operatorname{Spin}^{\C}$) structure, (see the discussion before Theorem~\ref{main}).  
In the next example we give the list of all 5-dimensional real Bott manifolds (with) without $\operatorname{Spin} (\operatorname{Spin}^{\C}$) structure.
\begin{ex}\label{ex2}
From {\em \cite{Nazra}} we have the list of all 5-dimensional oriented real Bott manifolds. There are 7 such manifolds without the torus.
%Moreover the oriented real Bott manifold $M(C),$ where
Here are their matrices:
$$A_4 =  \left[
\begin{matrix}
0&1&0&1&0\\
0&0&1&0&1\\
0&0&0&1&1\\
0&0&0&0&0\\
0&0&0&0&0
\end{matrix}\right], 
A_{23} = \left[
\begin{matrix}
0&1&1&0&0\\
0&0&0&0&0\\
0&0&0&1&1\\
0&0&0&0&0\\
0&0&0&0&0
\end{matrix}\right],$$

$$A_{29} = \left[
\begin{matrix}
0&1&1&1&1\\
0&0&0&0&0\\
0&0&0&0&0\\
0&0&0&0&0\\
0&0&0&0&0
\end{matrix}\right],
A_{37} = \left[
\begin{matrix}
0&0&0&0&0\\
0&0&1&1&0\\
0&0&0&1&1\\
0&0&0&0&0\\
0&0&0&0&0
\end{matrix}\right],$$

$$A_{40} = \left[
\begin{matrix}
0&0&1&0&1\\
0&0&1&1&0\\
0&0&0&1&1\\
0&0&0&0&0\\
0&0&0&0&0
\end{matrix}\right],
A_{48} = \left[
\begin{matrix}
0&0&1&0&1\\
0&0&1&1&0\\
0&0&0&0&0\\
0&0&0&0&0\\
0&0&0&0&0
\end{matrix}\right],$$

$$A_{49} = \left[
\begin{matrix}
0&0&0&0&0\\
0&0&0&0&0\\
0&0&0&1&1\\
0&0&0&0&0\\
0&0&0&0&0
\end{matrix}\right].$$
From the first part of Theorem~\ref{main} above, for $i = 1, j = 2$ the manifold $M(A_{4})$ 
has no $\operatorname{Spin^{\C}}$ structure.
For the same reasons (for $i = 1, j = 2$) manifolds $M(A_{40})$ and $M(A_{48})$ 
have no $\operatorname{Spin}$ structures.
The manifold $M(A_{23})$ has no a $\operatorname{Spin}$ structure, because it satisfies 
for $i = 1, j = 3$ the second part of the {\em Theorem~\ref{main}}.
Since any flat oriented manifold with $\Z_2$ holonomy
has $\operatorname{Spin}$ strucure, {\em \cite[Theorem 3.1]{HS}} manifolds $M(A_{29}), M(A_{49})$ have it. Our last example, the manifold $M(A_{37})$ has $\operatorname{Spin}$   
structure and we leave it as an exercise.
\vskip 2mm
\noindent
In all these cases it is possible to calculate the $w_{2}$ with the help 
of (\ref{w1}), (\ref{sw1}) and (\ref{ring}). In fact,
$w_2(M(A_4)) = (x_2)^2 + x_{1}x_{3}, w_2(M(A_{23})) = x_{1}x_{3}, w_2(M(A_{40})) = w_2(M(A_{48})) = x_{1}x_{2}.$ In all other cases
$w_2 = 0.$
\end{ex}
\begin{ex}
Let 
$$A = \left[
\begin{matrix}
0&0&1&1&1&1&0\\
0&0&0&0&0&1&1\\
0&0&0&\ast&\ast&\ast&\ast\\
0&0&0&0&\ast&\ast&\ast\\
0&0&0&0&0&\ast&\ast\\
0&0&0&0&0&0&0\\
0&0&0&0&0&0&0
\end{matrix}\right],$$
be a family of Bott matrices, with $\ast\in\Z_2.$
It is easy to check that the first two rows
satisfy the condition of Theorem \ref{main}. Hence the oriented real 
Bott manifolds $M(A)$ have no the $\operatorname{Spin}$ structure. 
\vskip 1mm
\noindent 
\end{ex}
\begin{rem}\label{remark}
In {\em \cite{AS}} on page 6 an example of the flat (real Bott) manifold 
$M$ without $\operatorname{Spin}$ structure is considered.
By an immediate calculation the Bott matrix of $M$ is equal to
$$\left[
\begin{matrix}
0&0&1&1&0\\
0&0&0&1&1\\
0&0&0&0&0\\
0&0&0&0&0\\
0&0&0&0&0
\end{matrix}\right].$$   
\end{rem}
\section{Concluding Remarks}
The tower (\ref{tower}) is an analogy of a Bott tower
$$
W_n\to W_{n-1}\to ...\to W_1=\C P^1\to W_0 =\{\text{a}\hskip 1mm\text{point}\}
$$ where $W_i$ is a $\C P^1$ bundle on $W_{i-1}$
i.e.; $W_{i} = P(1\oplus L_{i-1})$ and $L_{i-1}$ is a holomorphic
line bundle over $W_{i-1}.$ As in (\ref{tower}) $P(1\oplus L_{i-1})$ 
is projectivisation of the trivial linear bundle and $L_{i-1}$.
It was introduced by Grossberg and Karshon \cite{GK}.
As is well known, see \cite{CMO} for the complete bibliography, $W_n$ is a toric manifold.
\vskip 5mm
\noindent
There is an open problem: Is it true that two toric manifolds are diffeomorphic (or homeomorphic) if their cohomology
rings with integer coefficients are isomorphic as graded rings ?
In some cases it has partial affirmative solutions (see \cite{KM}).
For real Bott manifolds the following is true.
\vskip 1mm
\noindent
{\bf Theorem} (\cite[Theorem 1.1]{KM}) {\em Two real Bott manifolds are diffeomorphic if 
and only if their cohomology rings with
$\Z_2$ coefficients are isomorphic as graded rings. Equivalently, they are cohomological rigid.}
\vskip 1mm
\noindent
All of this suggests the following:
\vskip 2mm
\noindent
{\bf Question} Are ${\cal GHW}$-manifolds cohomological rigid ?
\vskip 2mm
\noindent
The answer to the above question is positive for manifolds from ${\cal GHW}\cap {\cal RBM}$.
It looks the most interesting for oriented GHW-manifolds. However,
for $n = 5$ there are two oriented Hantzsche-Wendt manifolds. From direct calculations with the help of a computer we know that they have different
cohomology rings with $\Z_2$ coefficients.

\vskip 2mm
\noindent
Institute of Mathematics, Maria Curie-Sk{\l}odowska University\\
Pl. Marii Curie-Sk{\l}odowskiej 1\\
20-031 Lublin\\
Poland\\
E-mail: anna.gasior @ poczta.umcs.lublin.pl
\vskip 5mm
\noindent
Institute of Mathematics, University of Gda\'nsk\\
ul. Wita Stwosza 57,\\
80-952 Gda\'nsk,\\
Poland\\
E-mail: matas @ univ.gda.pl
\end{document}